\documentclass{LMCS}

\def\dOi{10(4:7)2014}
\lmcsheading%
{\dOi}
{1--9}
{}
{}
{Feb.~13, 2013}
{Dec.~\phantom09, 2014}
{}

\ACMCCS{[{\bf Mathematics of computing}]: Continuous
  mathematics---Continuous functions; [{\bf Theory of computation}]:
  Design and analysis of algorithms---Approximation algorithms
  analysis---Numeric approximation algorithms; Computational
  complexity and cryptography}


\usepackage{hyperref}

\theoremstyle{definition}

\newtheorem{remark}{Remark}[section]

\newcommand{\In}{\subseteq}

\newcommand{\IN}{{\mathbb N}}

\newcommand{\IR}{{\mathbb R}}
\newcommand{\IC}{{\mathbb C}}

\begin{document}

\title[Computing a Solution of Feigenbaum's Functional Equation in
  Polynomial Time] {Computing a Solution of Feigenbaum's Functional
  Equation in Polynomial Time}

\author[P.~Hertling]{Peter Hertling}
\address{
Computer Science Department,
Universit\"at der Bundeswehr M\"unchen,
85577 Neubiberg,
Germany}
\email{\{peter.hertling,christoph.spandl\}@unibw.de}

\author[C.~Spandl]{Christoph Spandl}
\address{}
\email{}


\begin{abstract}
Lanford has shown in 1982 that Feigenbaum's functional equation has an
analytic solution.  We show that this solution is a polynomial time
computable function.  This implies in particular that the so-called
first Feigenbaum constant is a polynomial time computable real number.
\end{abstract}

\maketitle

\section{Introduction}

Independently, Feigenbaum~\cite{Fei78} and Gro{\ss}mann and Thomae~\cite{GT77} observed that the
behaviour of the points of bifurcations of certain parameterized classes of dynamical systems on an interval
obeys certain universal laws that are governed by constants which are now called Feigenbaum constants.
For detailed presentations of these notions the reader is referred to \cite{CE80} and to \cite{CEL80}.
In particular the so-called first Feigenbaum constant $\alpha = -2.50290787\ldots$
is the inverse $1/g(1)$ of the value $g(1)$ at $1$ of a solution $g$ of Feigenbaum's functional equation
which was explicitly constructed by Lanford~\cite{Lan82}.
In this note we show that this solution function $g$ is a polynomial time computable function.
This implies that the first Feigenbaum constant is a polynomial time computable number.

Which real numbers are computable? This question was one of the motivations for Alan Turing
to write his famous papers \cite{Tur36,Tur37}, in which he developed the notion of a Turing machine
and gave a definition of computable real numbers. Since then computable analysis has developed into a research area
in which the effective solvability of problems over the real numbers or more general continuous objects,
in particular all kinds of numerical problems, is analyzed using mathematically precise notions of effective
solvability, based on computability theory and complexity theory; see, e.g., \cite{PR89,Ko91,Wei00,BHW08}.
Among the first questions that one can ask in this theory is the question whether particular real number
constants are computable real numbers or not. For example, it is easy to see and well known
that the number $\pi$ and the Euler number $\mathrm{e}$ are computable.
In fact, they can be computed quite fast. An exemplary recent result of this kind is the observation
by Rettinger~\cite{Ret08} that the Bloch constant, a famous real number constant in complex analysis, is
computable.

A real number $c$ is called {\em computable}, if
there is an algorithm (a Turing machine) which, given an arbitrary $n\in\IN$
computes a rational number $q_n$ satisfying $|c - q_n| < 2^{-n}$.
A real number $c$ is called {\em polynomial time computable} if there are a Turing machine $M$ and a
polynomial $p$ with coefficients in $\IN$ such that $M$, given the string $1^n$ for any $n\in\IN$,
computes in at most $p(n)$ steps a binary string $a=a_m\ldots a_0$ (where $m$ is an arbitrary natural number)
and a binary string $b = b_1\ldots b_n$ such that
\[ |c - a.b| < 2^{-n} .\]
Here, by $a.b$ we mean the dyadic rational number defined by
\[ a.b
  = a_m\ldots a_0 . b_1\ldots b_n
  = \sum_{i=0}^m a_i \cdot 2^i + \sum_{j=1}^n b_j \cdot 2^{-j} . \]
Instead of binary strings one might as well consider decimal strings,
and instead of the upper bound $2^{-n}$ one might as well consider $10^{-n}$.
Finally, a sequence $(c_k)_{k\in\IN}$ of real numbers is called {\em polynomial time computable}
if there are a Turing machine $M$ and a two-variate
polynomial $p$ with coefficients in $\IN$ such that $M$, given $1^k01^n$ for $k,n\in\IN$,
computes in at most $p(k,n)$ steps a binary string $a=a_m\ldots a_0$ (where $m$ is an arbitrary natural number)
and a binary string $b = b_1\ldots b_n$ such that
\[ |c_k - a.b| < 2^{-n} .\]

In order to formulate our main result precisely we need to introduce some terminology.
We closely follow Lanford~\cite{Lan82}. In fact, this paper by Lanford is the basis of our analysis.

Let $M$ be the set of all continuously differentiable functions $f:[-1,1]\to[-1,1]$ satisfying the following conditions:
\begin{enumerate}
\item
$f(0)=1$,
\item
$x \cdot f^\prime(x) < 0$ for $x \neq 0$, i.e., $f$ is strictly increasing on $[-1,0]$ and strictly decreasing on $[0,1]$,
\item
$f(-x)=f(x)$ for all $x$, i.e., $f$ is even.
\end{enumerate}
Furthermore, let $D \In M$ be the set of all functions in $M$ satisfying additionally the following conditions:
\begin{enumerate}
\item
$0 < -f(1)$,
\item
$-f(1) < f(f(1))$,
\item
$f(f(f(1))) \leq -f(1)$.
\end{enumerate}
It is easy to check that for any function $f\in D$, the function $Tf$, defined by
\[ Tf(x) := \frac{1}{f(1)} \cdot f( f( f(1) \cdot x) ) \]
is an element of $M$.
Lanford~\cite{Lan82} showed the following result.

\begin{thm}[{\cite[Theorem 1 and Prop.~2]{Lan82}}]
\label{theorem:1-lanford}
There exists a function $g$, analytic and even on the set $\{z \in \IC ~:~ |z| < \sqrt{8}\}$
and with real values on real numbers, whose restriction to $[-1,1]$ is an element of $D$ and
a fixed point of the operator $T$.
\end{thm}

The so-called first Feigenbaum constant $\alpha$ is given by $\alpha = 1/g(1)$.

We prove the following addition to Lanford's theorem.

\begin{thm}
\label{theorem:main}
There exists a function $g$ that has the properties stated in Theorem~\ref{theorem:1-lanford}
and additionally the following properties.
\begin{enumerate}\enlargethispage{4\baselineskip}
\item
The sequence of Taylor coefficients around $0$ of this analytic function $g$ is a
polynomial time computable sequence of real numbers.
\item
The number $\alpha = 1/g(1)$ is a polynomial time computable real number.
\end{enumerate}
\end{thm}

Our proof is based on Lanford's paper~\cite{Lan82}.
In the following section we give the proof.

\section{A Polynomial Time Algorithm for Computing Lanford's solution of Feigenbaum's Functional Equation}

Lanford uses a variant of the Newton method in order to define an operator which has the same
fixed points as $T$. Then he gives a computer-assisted proof of a number of estimates
that show that this operator is a contraction in the neighborhood of an explicitly defined polynomial $\psi_0$,
with respect to an $\ell_1$-type norm on the space of Taylor coefficients of functions closely related to the functions
$f$ on which $T$ acts.
Furthermore, this operator maps this polynomial $\psi_0$ not too far away from itself.
By the contraction mapping principle it follows that the operator has a unique fixed point $g$.

We show that this construction leads to a polynomial time algorithm.

The following terminology is copied from \cite{Lan82}.
Let
\[ \Omega:= \{z \in \IC ~:~ |z^2-1| < 2.5\} , \]
and let $\mathcal{H}$ be the Banach space of even functions, bounded and analytic on $\Omega$,
real on real points, equipped with the supremum norm. We also define
\begin{eqnarray*}
  \mathcal{H}_0 &:=& \{ f \in \mathcal{H} ~:~ f(0)=0\}, \\
  \mathcal{H}_1 &:=& 1 + \mathcal{H}_0 .
\end{eqnarray*}
Lanford works on a subspace of $\mathcal{H}_1$ equipped with a stronger norm.
Let $\IN^+ := \{1,2,3,\ldots\}$ be the set of positive integers,
$\ell_1 := \{ \nu:\IN^+ \to \IR ~:~ \sum_{i=1}^\infty |\nu_i| < \infty\}$,
where $\nu_i := \nu(i)$.
Lanford considers the space $\IR\oplus \ell_1$ with the $\ell_1$-norm $||\cdot ||$ defined by
\[ ||(u,\nu)|| := |u| + \sum_{i=1}^\infty |\nu_i|  \]
for $u\in\IR$ and $\nu\in\ell_1$.
With an element $(u,\nu)\in \IR\oplus \ell_1$ he associates
the following element $\psi$ of $\mathcal{H}_1$:
\[ \psi(z) := 1 - z^2 \cdot \left( \frac{u}{10}
              + \sum_{i=1}^\infty \nu_i \cdot \left(\frac{z^2-1}{2.5}\right)^i \right) .  \]
Let $\mathcal{A}$ be the set of all functions $\psi$ defined in this way.
$\mathcal{A}$ is a subset of $\mathcal{H}_1$ and contains any element of $\mathcal{H}_1$
that is analytic on a neighborhood of the closure of $\Omega$.
In the following we will identify the elements of $\mathcal{A}$ with elements of
the space $\IR\oplus \ell_1$ with the $\ell_1$-norm introduced above.

The first step in Lanford's construction is the explicit definition of a polynomial
$\psi_0 \in \mathcal{A}$ of degree $20$ of the form
$\psi_0(z) = 1 + \sum_{i=1}^{10} g_i^{(0)} \cdot z^{2i}$
by choosing as the values $(g_1^{(0)},\ldots,a_{10}^{(0)})$:
``the first ten terms of the series given in Table 1 below'';
this table can be found on Page 432 in \cite{Lan82}.
Then Lanford continues by stating that for $\psi\in\mathcal{A}$ with
$||\psi - \psi_0|| < 0.01$ one has $T\psi \in \mathcal{A}$ as well.
The goal is to compute a fixed point of $T$ as the limit of a sequence of functions
starting with $\psi_0$ that are computed using a contractive mapping.
In order to achieve this, Lanford uses the operator $J:\IR\oplus \ell_1\to\IR\oplus \ell_1$
defined by
\[ J(u,\nu) := \left( \frac{u}{3.669}, - \nu \right) , \]
and defines for any $\psi \in \mathcal{A}$ with $||\psi - \psi_0||<0.01$
\[ \Phi(\psi) := \psi - J( T\psi - \psi) . \]
This operator $\Phi$ is an approximation of the operation iterated
in the Newton algorithm applied to the
function $\psi \mapsto T\psi - \psi$.
Note that $\Phi$ has the same fixed points as $T$.
For the proof of the following estimates Lanford uses computer calculations.
By $D\Phi(\psi)$ in the following lemma we mean the Fr\'{e}chet derivative of
$\Phi$ at $\psi$, which exists and can easily be calculated.

\begin{lem}[{\cite[Page 430]{Lan82}}]
\label{lemma:calc}\hfill
\begin{enumerate}
\item
For $\psi \in\mathcal{A}$ with $||\psi - \psi_0|| < 0.01$,
$|| D\Phi(\psi) || < 0.9$.
\item
$||\Phi(\psi_0) - \psi_0|| < 4 \cdot 10^{-6}$.
\end{enumerate}
\end{lem}


This lemma implies that $\Phi$ maps the closed ball
$\{\psi \in \mathcal{A} ~:~ ||\psi - \psi_0|| \leq 0.009\}$ into itself and that $\Phi$ is a 
contraction with Lipschitz constant $0.9$ on this ball.
By the contraction mapping theorem, the sequence $(\phi_m)_m$ of functions defined
by $\phi_0:=\psi_0$ and
$\phi_{m+1} := \Phi(\phi_m)$
converges to a fixed point $g$ of $\Phi$. It satisfies
\begin{equation}
\label{equation:contraction}
 || g - \phi_m|| \leq 4 \cdot 10^{-5} \cdot 0.9^m , 
\end{equation}
hence,
\begin{equation}
\label{equation:contraction0}
 ||g - \psi_0|| \leq 4 \cdot 10^{-5} .
\end{equation}
Remember that $\Phi$ has the same fixed points as $T$.
Thus, $g$ is a fixed point of $T$.
Lanford shows that this fixed point of $T$ has all of the properties
stated in Theorem~\ref{theorem:1-lanford}.

{}From~\eqref{equation:contraction}
it is clear that by starting with the explicitly defined polynomial $\psi_0$ and by
applying the contractive operator $\Phi$ to $\psi_0$ \ $O(n)$ times
one obtains a polynomial that approximates the fixed point $g$
with precision $10^{-n}$ (with respect to the norm considered by Lanford and described
above).
We wish to show that one can approximate the fixed point $g$
with precision $10^{-n}$ in time polynomial in $n$.
In order to achieve that, we are going to show that the precision needed in the $n$-th step is not too high
and that the number of coefficients that need to be considered in the $n$-th step is not too high as well.
In fact, we will show that in the $n$-th step it is sufficient to consider a polynomial of a degree depending linearly on $n$. 

First, we make some observations about the fixed point $g$ of the operators $T$ and $\Phi$.
For $z_0 \in \IC$ and $r>0$ let
\[ B(z_0,r) := \{ z \in \IC ~:~ |z-z_0| < r\}. \]
By Theorem~\ref{theorem:1-lanford}, $g$ is an even analytic function
defined on the disc $B(0,\sqrt{8})$ satisfying $g(0)=1$.
Therefore its Taylor series around $0$,
\[ 1 + \sum_{i=1}^\infty a_i z^{2i} \text{ with } a_i := \frac{g^{(2i)}(0)}{(2i)!} \]
converges in $B(0,\sqrt{8})$ and is equal to $g$ in $B(0,\sqrt{8})$.
Then the function $h$ defined by
\[ h(z) := - \sum_{j=0}^\infty a_{j+1} z^j \]
is an analytic function in the ball $B(0,8)$, and for all $z\in B(0,\sqrt{8})$
we have
\[ g(z) = 1 - z^2 \cdot h(z^2) . \]
The Taylor series of $h$ around $1$ converges and is identical with $h$ in the ball $B(1,7)$:
\begin{equation*}
h(z) = \sum_{j=0}^\infty b_j(z-1)^j \text{ with } b_j:= \frac{h^{(j)}(1)}{j!}.
\end{equation*}
The Cauchy integral formula then gives
\begin{equation*}
|b_j|\leq C\cdot 6.5^{-j}
\end{equation*}
with $C=\max\{|h(z)| : z\in\partial B(1,6.5)\}$,
where for $z_0\in \IC$ and $r>0$
\[ \partial B(z_0,r) = \{z \in \IC ~:~ |z-z_0|=r\} . \]
We claim that results in \cite{Lan82} imply $C \leq 62/13$. 
Indeed, according to \cite[Remark 4.2 on Page 431]{Lan82}
$|g(z)-g^{(0)}(z)|\leq 1.7\cdot 10^{-2}$ holds for all
$z\in B(0,\sqrt{8})$, where $g^{(0)}(z)=1+\sum_{i=1}^{40}g^{(0)}_i\cdot z^{2i}$
is a polynomial of degree $80$ with coefficients $g^{(0)}_i$ given in Table 1
on Page 432 of \cite{Lan82}. Defining $h^{(0)}(z)$ by $g^{(0)}(z)=1-z^2\cdot h^{(0)}(z^2)$,
we obtain
$|h(z)-h^{(0)}(z)|\leq 1.7\cdot 10^{-2} / 5.5  \leq 3.1 \cdot 10^{-3}$ 
for $z\in\partial B(1,6.5)$.
Writing
$h^{(0)}(z) = \sum_{j=0}^{39} b_j^{(0)}(z-1)^j$
(the numbers $b_j^{(0)}$ for $j=0,\ldots,39$ 
can easily be computed from the numbers $g_i^{(0)}$ for $i=1,\ldots, 40$)
we obtain
\begin{equation*}
|h(z)|\leq 3.1\cdot 10^{-3} + |h^{(0)}(z)|\leq
3.1\cdot 10^{-3} + \sum_{j=0}^{39}|b^{(0)}_j|\cdot 6.5^j 
\leq 4.7 < \frac{62}{13}
\end{equation*}
for $z\in\partial B(1,6.5)$. Thus, $C \leq 62/13$, and we have
\begin{equation*}
|b_j|\leq \frac{62}{13}\cdot 6.5^{-j} .
\end{equation*}
Defining
\[ u^{(\infty)} := 10 \cdot b_0 \text{ and, for } i\in\IN^+, \
   \nu_i^{(\infty)} := 2.5^i \cdot b_i \]
we obtain for all $z\in\{ z \in \IC ~:~ |z^2 - 1| < 7\}$
\[ g(z) = 1 - z^2 \cdot \left( \frac{u^{(\infty)}}{10} +
   \sum_{i=1}^\infty \nu_i^{(\infty)} \cdot \left(\frac{z^2-1}{2.5} \right)^i \right)  \]
and for all $i\geq 1$
\begin{equation}
\label{equation:g-coeff-1}
 |\nu_i^{(\infty)}| \leq \frac{62}{13} \cdot \left( \frac{2.5}{6.5} \right)^i
     = \frac{62}{13} \cdot \left( \frac{5}{13} \right)^i .
\end{equation}
Note that this implies for all $k\geq 1$
\begin{equation}
\label{equation:g-coeff-2}
  \sum_{i=k}^\infty |\nu_i^{(\infty)}| \leq
   \frac{62}{13} \cdot \frac{13}{8} \cdot \left( \frac{5}{13} \right)^k =
   \frac{31}{4} \cdot \left( \frac{5}{13} \right)^k
\end{equation}

We wish to approximate in time polynomial in $n$
the function $g$, i.e., the sequence
$(u^{(\infty)},\nu^{(\infty)}) \in \IR\oplus\ell_1$, with precision $10^{-n}$ in the norm introduced
above. We start with the polynomial $\psi_0$ chosen by Lanford and define the numbers
$u^{(0)},\nu_1^{(0)},\ldots,\nu_{9}^{(0)}$ by
\[ \psi_0(z) := 1 - z^2 \cdot \left( \frac{u^{(0)}}{10}
              + \sum_{i=1}^{9} \nu_i^{(0)} \cdot \left(\frac{z^2-1}{2.5}\right)^i \right).   \]
These numbers can easily be computed explicitly
and are given in Table~\ref{table:coeff}.

\[
\begin{array}{c|r}
\multicolumn{2}{c}{\text{Table~\ref{table:coeff}}} \\
\hline
u^{(0)} &   13.99535280247654509657069657886239000000000 \\
\nu_1^{(0)} &  - 0.37020336425570944099807863650264000000000 \\
\nu_2^{(0)} &  - 0.10516441308487059395306704671240000000000 \\ 
\nu_3^{(0)} &    0.04689224531866417356902064258837500000000 \\
\nu_4^{(0)} &  - 0.00657196434429489515940234119726562500000 \\
\nu_5^{(0)} &  - 0.00092424880356949042888086870078125000000 \\
\nu_6^{(0)} &    0.00060199775715465703408272872656250000000 \\
\nu_7^{(0)} &  - 0.00007266358160903580114416214843750000000 \\
\nu_8^{(0)} &  - 0.00003921160572782132082950382843017578125 \\
\nu_9^{(0)} &    0.00000105783506805382222151565551757812500
\end{array}
\label{table:coeff}
\]\medskip

We come to the \\
{\bf Central Step} of the algorithm.
Let us assume by induction hypothesis that for some $m\geq 0$ we have computed
$10+m$
real numbers $u^{(m)},\nu_1^{(m)},\ldots,\nu_{9+m}^{(m)}$
with the following properties:
\begin{itemize}
\item
{\bf Property I:}
Each of these numbers
is a decimal fraction of the form
\[ \sigma a_1 a_0 . b_1 \ldots b_{41+m} \]
where $\sigma \in \{-1,+1\}$ is a sign and
$a_1,a_0,b_1,\ldots, b_{41+m}$ are decimal digits.
\item
{\bf Property II:}
The polynomial $\psi_m$ of degree $20 + 2m$ defined by
\[ \psi_m(z) := 1 - z^2 \cdot \left( \frac{u^{(m)}}{10}
              + \sum_{i=1}^{9+m} \nu_i^{(m)} \cdot \left(\frac{z^2-1}{2.5}\right)^i \right)   \]
satisfies
\[ ||\psi_m - g|| < 0.01 \cdot 0.93^m . \]
\end{itemize}
First we observe that the polynomial $\psi_0$ indeed has these properties for $m=0$.
The first property
can be checked easily by explicitly calculating the numbers
$u^{(0)},\nu_1^{(0)},\ldots,\nu_9^{(0)}$; see Table~\ref{table:coeff}.
The second property is a consequence of Equation~\eqref{equation:contraction0}.

Now we wish to compute suitable numbers
$u^{(m+1)},\nu_1^{(m+1)},\ldots,\nu_{9+m+1}^{(m+1)}$.
Note that, due to \eqref{equation:contraction0}, Property II implies
\begin{equation}
\label{equation:psin-psi0}
 ||\psi_m - \psi_0||
 \leq ||\psi_m - g|| + ||g - \psi_0||
 < 0.01 \cdot 0.93^m + 4 \cdot 10^{-5} < 0.01
\end{equation}
if $m\geq 1$.
Therefore, the estimate stated in Lemma~\ref{lemma:calc}(1) applies to $\psi_m$.
As $\psi_m$ is a polynomial of degree $\leq 20 + 2m$ in which no
monomials of odd degree occur,
by definition of $\Phi$ and of $T$,
$\Phi(\psi_m)$ is a polynomial of degree $\leq (20 + 2m)^2$
in which no monomials of odd degree occur.
The numbers $v,\mu_1,\ldots,\mu_{2\cdot (10+m)^2-1}$ are defined by
\[ \Phi(\psi_m)(z) = 1 - z^2 \cdot \left( \frac{v}{10}
              + \sum_{i=1}^{2\cdot (10+m)^2-1} \mu_i \cdot \left(\frac{z^2-1}{2.5}\right)^i \right).   \]
In order to simplify notation in the following, we additionally use
\[ \mu_i := 0 \text{ for } i > 2\cdot (10+m)^2-1  . \]
We compute the first $10 + m+1$ of these numbers with precision $10^{-41-(m+1)}$, i.e.,
we compute decimal fractions 
$u^{(m+1)},\nu_1^{(m+1)},\ldots,\nu_{9+m+1}^{(m+1)}$
with at most $41+(m+1)$ digits after the decimal point such that
\[ |u^{(m+1)} - v| \leq 10^{-41-(m+1)} \text{ and for }
  i=1,\ldots,9+m+1, \quad
  |\nu_i^{(m+1)} - \mu_{i}| \leq 10^{-41-(m+1)} . \]
That is, we simply forget the coefficients $\mu_i$ for $i > 9+m+1$.
This ends the description of the Central Step of the algorithm in which we compute
the numbers $u^{(m+1)},\nu_1^{(m+1)},\ldots,\nu_{9+m+1}^{(m+1)}$ from the numbers
$u^{(m)},\nu_1^{(m)},\ldots,\nu_{9+m}^{(m)}$.

It is clear that this Central Step can be executed in time polynomial in $m$. 

\begin{remark}
\label{remark:precise-complexity}
It is fairly easy to see that no more than $O(m^3)$ elementary arithmetic
operations are needed in the Central Step. Note that in order to achieve that it is important
that not all of the $2\cdot (10+m)^2$ coefficients of the polynomial
$\Phi(\psi_m)(z)$ are computed but only the first $10+m+1$ coefficients.
By somewhat tedious estimations one can show that there
are positive constants $a,b$ with the property that it is sufficient to 
perform each arithmetic operation with $a+b\cdot m$ digits in total, that is,
before or after the decimal point.
Let $\mathcal{M}(m)$ be a function satisfying
$O(\mathcal{M}(c\cdot m)) = O(\mathcal{M}(m))$
such that two binary or decimal
numbers of length $m$ can be multiplied in time $O(\mathcal{M}(m))$.
For example, the Sch{\"o}nhage-Strassen bound $m \cdot \log m \cdot \log\log m$ is such a function.
It is well known that one can also add, subtract, or divide numbers of length $m$ within this time~\cite{Bre76}.
We conclude that the Central Step can be done in time $O(m^3 \cdot \mathcal{M}(m))$.
\end{remark}


We claim that these numbers $u^{(m+1)},\nu_1^{(m+1)},\ldots,\nu_{9+m+1}^{(m+1)}$
have the desired properties.
First, we show that they have Property II (with $m$ replaced by $m+1$).
We calculate \enlargethispage{\baselineskip}
\begin{align*}
 \lefteqn{||\psi_{m+1} - g||} && \\
 & = |u^{(m+1)} - u^{(\infty)}|
    + \sum_{i=1}^{9+m+1} |\nu_i^{(m+1)} - \nu_i^{(\infty)}|
    + \sum_{i=10+m+1}^\infty |\nu_i^{(\infty)}| \\
 & \leq |u^{(m+1)} - v| + |v - u^{(\infty)}|
    + \sum_{i=1}^{9+m+1} |\nu_i^{(m+1)} - \mu_i|
    + \sum_{i=1}^{9+m+1} |\mu_i - \nu_i^{(\infty)}|
    + \sum_{i=10+m+1}^\infty |\nu_i^{(\infty)}| \\
%
 & \leq (10 + m+1) \cdot 10^{-41-(m+1)}
    + |v - u^{(\infty)}|
    + \sum_{i=1}^{\infty} |\mu_i - \nu_i^{(\infty)}|
    + \sum_{i=10+m+1}^\infty |\nu_i^{(\infty)}| \\
 & = (10 + m+1) \cdot 10^{-41-(m+1)}
    + ||\Phi(\psi_m) - g|| + \sum_{i=10+m+1}^{\infty} |\nu_i^{(\infty)}| \\   
 & = (10 + m+1) \cdot 10^{-41-(m+1)}
    + ||\Phi(\psi_m) - \Phi(g)|| + \sum_{i=10+m+1}^{\infty} |\nu_i^{(\infty)}| \\
 & \leq (10 + m+1) \cdot 10^{-41-(m+1)} + 0.9 \cdot ||\psi_m - g||
     + \frac{31}{4}\cdot \left( \frac{5}{13} \right)^{10+m+1} \\
 & < (10 + m+1) \cdot 10^{-41-(m+1)} + 0.9 \cdot 0.01 \cdot 0.93^m
     + \frac{31}{4}\cdot \left( \frac{5}{13} \right)^{10+m+1} \\
 & \leq 0.01\cdot0.001\cdot0.93^m+0.01\cdot0.9\cdot0.93^m+0.01\cdot0.025\cdot0.93^m\\
 & < 0.01 \cdot 0.93^{m+1}.
\end{align*}
The next to last estimate is a consequence of the bounds
\begin{align*}
(10 + m+1) \cdot 10^{-41-(m+1)} &\leq  10^{-5} \cdot 0.93^{m} \text{ and } \\
\frac{31}{4} \cdot \left( \frac{5}{13} \right)^{10+m+1} &\leq
  2.5 \cdot 10^{-4} \cdot 0.93^{m}
\end{align*}
for $m\geq 0$. It is not hard to verify them. Thus we have shown
Property II (with $m$ replaced by $m+1$ in Property II) for the numbers
$u^{(m+1)},\nu_1^{(m+1)},\ldots,\nu_{9+m+1}^{(m+1)}$.
Concerning Property I, it is part of the definition of these numbers,
that each of them is a finite decimal fraction with at most $41+(m+1)$ decimal digits
after the decimal point. That each of them has at most $2$ decimal digits
in front of the decimal point follows from the values of the numbers
$u^{(0)},\nu_1^{(0)},\ldots,\nu_9^{(0)}$,
from Equation~\eqref{equation:contraction0} and Property II which together
imply that the numbers $u^{(m+1)},\nu_1^{(m+1)},\ldots,\nu_{9+m+1}^{(m+1)}$
are quite close to the numbers
$u^{(0)},\nu_1^{(0)},\ldots,\nu_9^{(0)},0, 0, 0, \ldots$;
see \eqref{equation:psin-psi0}.
We have shown that the numbers
$u^{(m+1)},\nu_1^{(m+1)},\ldots,\nu_{9+m+1}^{(m+1)}$
have Property I (with $m$ replaced by $m+1$ in Property I) as well.

Since the Central Step in the $m$-th iteration can be performed in time $O(m^3 \cdot \mathcal{M}(m))$,
the first $m$ iterations can be done in time $O(m^4 \cdot \mathcal{M}(m))$.
Since the smallest number $m$ with $0.01 \cdot 0.93^{m} \leq 10^{-n}$ depends linearly on $n$,
this shows that we can compute the first $n$ elements of the sequence $(u^{(\infty)},\nu^{(\infty)}) \in \IR\times \ell_1$
with precision $10^{-n}$ in time $O(n^4 \cdot \mathcal{M}(n))$.
Then also the number $g(1) = 1 - \frac{u^{(\infty)}}{10}$ can be computed in this time with precision $10^{-n}$.
The same applies to $1/g(1)$ because the inverse $1/y$ of a real number $y \neq 0$ can be computed with precision $10^{-n}$
in time $\mathcal{M}(n)$ using no more than a linear number of digits of $y$~\cite{Bre76,Wei00}.
Thus, the first Feigenbaum constant $\alpha=1/g(1)$ is a polynomial time computable
real number, actually, computable in time $O(n^4 \cdot \mathcal{M}(n))$.
%
%
By M{\"u}ller~\cite[Cor.~3.3]{Mue93} the fact that
the sequence $(u^{(\infty)},\nu^{(\infty)})$ can be computed in polynomial time
implies also that the function $g$ is computable
in polynomial time on any compact subset of the set
$\{ z \in \IC ~:~ |z^2 - 1| < 7\}$.
And by Theorem 2.4 of M{\"u}ller~\cite{Mue93}, the
sequence of Taylor coefficients around $0$ of the function $g$ is a
polynomial time computable sequence of real numbers.
We have shown Theorem~\ref{theorem:main}.

\section{Final Remarks}

The proof is based on a number of claims in Lanford's paper \cite{Lan82}.
In order to prove them, Lanford used analytic estimations and extensive computer calculations 
which are described only very roughly in his paper.
We must admit that we did not check all of these claims, but most of them.

Briggs~\cite{Briggs91,Bri97} has computed several hundred digits of the Feigenbaum constants,
but without correctness guarantee. We intend to
see how many digits of the first Feigenbaum constant
we can compute with a correctness guarantee,
using the exact real number arithmetic package
iRRAM by M{\"u}ller~\cite{Mue01}.


\end{document}